\documentclass[12pt]{article}

\newtheorem{theorem}{Theorem}[section]
\newtheorem{lemma}[theorem]{Lemma}
\newtheorem{proposition}[theorem]{Proposition}
\newtheorem{corollary}[theorem]{Corollary}

\newtheorem{definition}[theorem]{Definition}

\newenvironment{poof}{\textit{Proof:  }}{
~\hfill\rule{2mm}{3mm}\vspace{.2in}}

\usepackage{amssymb,amsmath,tabularx,graphicx,float,placeins}
\usepackage{tikz}
\usetikzlibrary[calc,intersections,through,backgrounds,arrows,decorations.pathmorphing]

\def\Acal{\mathcal{A}}

\def\Lcal{\mathcal{L}}

\def\Pcal{\mathcal{P}}

\def\Tcal{\mathcal{T}}

\def\Wcal{\mathcal{W}}

\def\1bb{\mathbb{1}}

\def\Rbb{\mathbb{R}}

\def\Zbb{\mathbb{Z}}

\def\ss{\scriptstyle}

\def\ra{\rightarrow}

\def\ov{\overline}

\def\wtil{\widetilde}

\def\pr{\prime}
\def\prpr{\prime\prime}

\def\setm{\setminus}

\DeclareMathOperator{\depth}{depth}
\DeclareMathOperator{\id}{id}

\begin{document}

\title{Biclosed Sets in Real Hyperplane Arrangements}
\date{\today}
\author{Thomas McConville}

\maketitle

\begin{abstract}
The set of chambers of a real hyperplane arrangement may be ordered by separation from some fixed chamber.  When this poset is a lattice, Bj\"orner, Edelman, and Ziegler proved that the chambers are in natural bijection with the biconvex sets of the arrangement.  Two families of examples of arrangements with a lattice of chambers are simplicial and supersolvable arrangements.  For these arrangements, we prove the that the chambers correspond to biclosed sets, a weakening of the biconvex property.
\end{abstract}

\section{Introduction}

The inversion set of a permutation of $\{1,\ldots,n\}$ is the collection of pairs $\{i,j\},\ 1\leq i<j\leq n$ such that $j$ precedes $i$ in the permutation.  Inversion sets of permutations are characterized as those collections $I$ of 2-element subsets of $\{1,\ldots,n\}$ for which
\begin{itemize}
\item if $\{i,j\}$ and $\{j,k\}$ are in $I$ then $\{i,k\}$ is in $I$, and
\item if $\{i,j\}$ and $\{j,k\}$ are not in $I$ then $\{i,k\}$ is not in $I$,
\end{itemize}
for $1\leq i<j<k\leq n$.

This characterization of inversion sets may be interpreted geometrically.  A finite central arrangement of hyperplanes in $\Rbb^n$ determines a complete fan of polyhedral cones, whose maximal cones are called chambers.  A permutation $\sigma$ of $\{1,\ldots,n\}$ corresponds to the chamber
$$c_{\sigma}=\{\mathbf{x}\in\Rbb^n:\ x_{\sigma(1)}<\cdots<x_{\sigma(n)}\}$$
of $\Acal_{braid}^n$, the arrangement of hyperplanes in $\Rbb^n$ defined by the equations $x_i=x_j$ for $i\neq j$.  The inversion set of a permutation $\sigma$ may be identified with the \emph{separation set} $S(c_{\sigma})$, the set of hyperplanes separating $c_{\sigma}$ from $c_{\id}$.  A subset $I$ of an arrangement $\Acal$ is \emph{biclosed} with respect to a chamber $c_0$ if $I\cap\Acal^{\pr}$ is the set of hyperplanes separating $c_0$ from some chamber of $\Acal^{\pr}$ for every rank 2 subarrangement $\Acal^{\pr}$.  By the characterization of inversion sets, a subset $I$ of $\Acal_{braid}^n$ is the separation set of some chamber if and only if $I$ is biclosed with respect to $c_{\id}$.  More generally, the elements of a finite Coxeter group $W$ correspond to biclosed subsets of its standard reflection arrangement $\Acal(W)$ (\cite{bourbaki:groupes}; see also the appendix of \cite{pilkington:convex}).  We prove this fact for a more general class of arrangements, which we recall in Sections \ref{sec_simplicial} and \ref{sec_supersolvable}.

\begin{theorem}\label{thm_main_biclosed}
Let $\Acal$ be a simplicial or supersolvable arrangement with fundamental chamber $c_0$.  For $I\subseteq\Acal$, $I$ is biclosed if and only if $I$ is the separation set of some chamber.
\end{theorem}

A hyperplane arrangement $\Acal$ with a fundamental chamber determines a convex geometry on $\Acal$, which we recall in Section \ref{sec_biclosed} (see \cite[Proposition 5.1]{bjorner.edelman.ziegler:lattice}).  The chambers of $\Acal$ form a poset $\Pcal(\Acal,c_0)$ where $c\leq d$ if $S(c)\subseteq S(d)$; see Figure \ref{fig_ex_chamber}.  If $\Acal$ is the standard reflection arrangement of a finite Coxeter group $W$, then this poset is isomorphic to the weak order on $W$.  Bj\"orner proved that the weak order of a finite Coxeter group is a lattice \cite{bjorner:orderings}.  Bj\"orner, Edelman and Ziegler proved that if $\Pcal(\Acal,c_0)$ is a lattice, then the chambers of $\Acal$ are in bijection with biconvex subsets of $\Acal$ \cite[Theorem 5.5]{bjorner.edelman.ziegler:lattice}.  Furthermore, for chambers $c$ and $d$ the separation set of $c\vee d$ is the convex closure of $S(c)\cup S(d)$.

\begin{figure}
\begin{centering}
\begin{tikzpicture}[scale=.5]
  \draw (-3,-4) -- (3,4)
        (-5,0) -- (5,0)
        (-3,4) -- (3,-4);
  \draw (0,-5/2) node[anchor=north]{$c_0$}
        (0,5/2) node[anchor=south]{$-c_0$};
  \begin{scope}[xshift=9cm]
    \draw (0,-5/2) node[anchor=north]{$c_0$} -- (-4/2,-3/2) -- (-4/2,3/2) -- (0,5/2) node[anchor=south]{$-c_0$} -- (4/2,3/2) -- (4/2,-3/2) -- cycle;
  \end{scope}
  \begin{scope}[xshift=16cm,yshift=-3cm,scale=2]
    \draw (0,0) -- (-1.5,1) -- (-1,2) -- (0,3) -- (1,2) -- (1.5,1) -- cycle;
    \draw (0,0) -- (0,1) -- (-1,2)
          (0,1) -- (1,2);
  \end{scope}
\end{tikzpicture}
\caption{\label{fig_ex_chamber}\scriptsize (left) An arrangement of three lines  (center) Poset of chambers (right) Poset of 2-closed subsets.}
\end{centering}
\end{figure}

A subset $I$ of an arrangement $\Acal$ with a fundamental chamber is \emph{2-closed} if $I\cap\Acal^{\pr}$ is convex for every rank 2 subarrangement $\Acal^{\pr}$ of $\Acal$.  By reduction to the rank 2 case, $I$ is biclosed if and only if both $I$ and $\Acal-I$ are 2-closed.  The 2-closure is typically weaker than convex closure and seldom defines a convex geometry, even for reflection arrangements \cite[Theorem 1(c)]{pilkington:convex}.  However, the two operators do agree for the reflection arrangements of types A or B \cite[Theorem 1(b)]{pilkington:convex}, \cite{stembridge:irreducible}, where the convex closure may be interpreted as a transitive closure for posets or signed posets, respectively (see \cite{reiner:signed}).  Stembridge determined the relative strength of various $r$-closures on reflection arrangements by computing their set of irreducible circuits, a special collection of circuits that suffices to compute the convex closure (see \cite[Proposition 1.1]{stembridge:irreducible} for a precise statement).

Even though the 2-closure is generally weaker than convex closure, Dyer proved for finite reflection arrangements that $S(c\vee d)$ is the 2-closure of $S(c)\cup S(d)$ for chambers $c$ and $d$ (\cite[Theorem 1.5]{dyer:weak}).  We prove an analogue of this result for certain hyperplane arrangements.  The bineighborly property will be recalled in Section \ref{sec_simplicial}.

\begin{theorem}\label{thm_main_2closure}
Let $\Acal$ be a bineighborly or supersolvable arrangement with fundamental chamber $c_0$.  If $c$ and $d$ are chambers of $\Acal$, then $S(c\vee d)$ is the 2-closure of $S(c)\cup S(d)$.
\end{theorem}

We would like to know whether the hypothesis of Theorem \ref{thm_main_2closure} may be replaced by a suitable lattice property of $\Pcal(\Acal,c_0)$.  In separate work, we proved that an arrangement $\Acal$ is bineighborly if and only if its chamber poset is a semidistributive lattice \cite[Theorem 4.2]{mcconville:crosscut}.  Semidistributivity was shown previously for the weak order of a finite Coxeter group by Le Conte de Poly-Barbut \cite{poly-barbut:treillis} and for chamber posets of simplicial arrangements by Reading \cite[Theorem 3]{reading:lattice}.

Reading proved that the poset of chambers of a supersolvable arrangement is a congruence-normal lattice \cite[Theorem 1]{reading:lattice}.  However, not all congruence-normal chamber lattices satisfy the conclusion of Theorem \ref{thm_main_2closure}.  For example, if $\Acal$ is a generic arrangement of four planes in $\Rbb^3$ with a simplicial fundamental chamber $c_0$, then $\Pcal(\Acal,c_0)$ is a congruence-normal lattice for which the 2-closure of the walls of $c_0$ is not the separation set of any chamber.

The set of reduced words of a Coxeter group element $w$ may be identified with maximal chains of the interval $[e,w]$ of the weak order.  A \emph{reduced gallery} between chambers $c_0$ and $c$ of an arrangement $\Acal$ is a maximal chain in the interval $[c_0,c]$ of $\Pcal(\Acal,c_0)$.  A maximal chain in $\Pcal(\Acal,c_0)$ induces a total order on $\Acal$.  We call a total order on $\Acal$ \emph{admissible} if its restriction to $\Acal^{\pr}$ defines a reduced gallery for all rank 2 subarrangements $\Acal^{\pr}$ of $\Acal$.

If $W$ is a (possibly infinite) Coxeter group, then an admissible total ordering of $\Acal(W)$ is called a \emph{reflection order} or \emph{convex order}.  When $W$ is finite, there is a well-known correspondence between reflection orders and reduced words for the longest element.  When $W$ is infinite, the collections of biclosed sets and reflections orders are not completely understood; see \cite{dyer:weak} or \cite{dyer.hohlweg.ripoll:imaginary} for conjectures and recent progress.  If $W$ is a Weyl group, a slightly different definition of convex order is used.  In this setting, convex orders for affine Weyl groups were classified by Ito \cite{ito:classification}.

The set of reduced galleries between a fixed pair of chambers forms a graph where two galleries are adjacent if one gallery may be obtained from the other by ``flipping'' about a codimension 2 intersection subspace.  The graph of reduced galleries was shown to be connected in successively greater generality by Tits \cite{tits:probleme}, Deligne \cite{deligne:immeubles}, Salvetti \cite{salvetti:topology}, and Cordovil-Moreira \cite{cordovil.moreira:homotopy}.  The graph of reduced galleries between opposite chambers has further topological connectivity (see \cite{bjorner:essential}) as well as further graph-theoretic connectivity (see \cite{athanasiadis.edelman.reiner:monotone} or \cite{athanasiadis.santos:monotone}).

More recently, the diameter of some reduced gallery graphs of supersolvable arrangements were computed by Reiner and Roichman \cite[Theorem 1.1]{reiner.roichman:diameter}.  Namely, if $c_0$ is incident to a modular flag of $\Acal$, then the graph of reduced galleries between $c_0$ and $-c_0$ has diameter equal to the number of codimension 2 intersection subspaces of $\Acal$.  In particular, the graph of reduced words for the longest element in types $A_n$ and $B_n$ have diameter in $O(n^4)$.  A key step in their proof relied on an unproven assumption, which we justify here by proving the following result.

\begin{theorem}\label{thm_main_gallery}
Let $\Acal$ be an arrangement with fundamental chamber $c_0$ for which every biclosed set is the separation set of some chamber.
\begin{enumerate}
\item A total order on $\Acal$ is admissible if and only if there exists a reduced gallery from $c_0$ to $-c_0$ inducing that order.
\item Let $r$ be a reduced gallery and $X$ is a codimension $2$ subspace.  If the total order defined by $r$ contains all of the hyperplanes of $\Acal_X$ consecutively, then there exists a reduced gallery with the same total order as $r$ with the hyperplanes in $\Acal_X$ reversed.
\end{enumerate}
\end{theorem}

The paper is structured as follows.  We recall some notation and fundamental results on real hyperplane arrangements in Section \ref{sec_chamber_posets}, following the notation of \cite{bjorner.lasVergnas.ea:oriented}.  We recall some facts about biclosed sets in Section \ref{sec_biclosed}.  We discuss reduced galleries and set-valued metrics on graphs in Section \ref{sec_galleries} and prove Theorem \ref{thm_main_gallery}.  In Section \ref{sec_simplicial}, we prove Theorem \ref{thm_main_biclosed} for simplicial arrangements and Theorem \ref{thm_main_2closure} for bineighborly arrangements.  We prove Theorems \ref{thm_main_biclosed} and \ref{thm_main_2closure} for supersolvable arrangements in Section \ref{sec_supersolvable}, and indicate the use of Theorem \ref{thm_main_gallery} in the diameter computation of the reduced gallery graph.

\section{Chamber Posets}\label{sec_chamber_posets}

A \emph{poset} is a set with a reflexive, symmetric, and transitive binary relation.  A \emph{lattice} is a poset for which every pair of elements $x,y$ has a least upper bound $x\vee y$ and a greatest lower bound $x\wedge y$.  The \emph{join} (\emph{meet}) of a finite subset $A$ of a lattice, denoted $\bigvee A$ ($\bigwedge A$), is the common least upper bound (greatest lower bound) of the elements of $A$.  If $x\leq y$, the \emph{closed interval} $[x,y]$ is the set of $z\in P$ with $x\leq z\leq y$.  If $x<y$, we say $y$ \emph{covers} $x$ if there does not exist $z$ for which $x<z<y$.

\begin{lemma}[\cite{bjorner.edelman.ziegler:lattice}, Lemma 2.1]\label{lem_local_lattice}
Let $P$ be a finite poset.  If $x\vee y$ exists whenever $x$ and $y$ both cover a common element, then $P$ is a lattice.
\end{lemma}

A \emph{real, central hyperplane arrangement} $\Acal$ is a finite set of (oriented) hyperplanes in $\Rbb^n$ whose common intersection contains the origin.  Each hyperplane $H$ partitions $\Rbb^n-H$ into two open half spaces, denoted $H^+$ and $H^-$.  For $x\in\{0,+,-\}^{\Acal}$, if $\bigcap_{H\in\Acal}H^{x(H)}$ is nonempty, then $\bigcap_{H\in\Acal}H^{x(H)}$ is a \emph{face} of $\Acal$ and $x$ is a \emph{covector} of $\Acal$.  We identify faces with their covectors when convenient.  For $x,y\in\{0,+,-\}^{\Acal}$, the \emph{composite} $x\circ y$ is the sign vector where for $H\in\Acal$
$$(x\circ y)(H)=\begin{cases}x(H)\ &\mbox{if}\ x(H)\neq 0\\y(H)\ &\mbox{if}\ x(H)=0\end{cases}.$$
The set of sign vectors $\{0,+,-\}^{\Acal}$ is given the product order where $0<+,\ 0<-$, and $+$ is incomparable with $-$.  For $x,y\in\{0,+,-\}^{\Acal}$, $y$ is \emph{incident} to $x$ if $x\leq y$.  A \emph{circuit} $v\in\{0,+,-\}^{\Acal}$ is a minimal sign vector such that $\bigcap_{\substack{H\in\Acal\\v(H)\neq 0}}H^{v(H)}$ is empty.

The set $\Lcal(\Acal)$ of covectors of an arrangement $\Acal$ satisfy

\begin{list}{}{}
\item[(L0)] $\mathbf{0}\in\Lcal$,
\item[(L1)] $x\in\Lcal$ implies $-x\in\Lcal$,
\item[(L2)] $x,y\in\Lcal$ implies $x\circ y\in\Lcal$, and
\item[(L3)] if $x,y\in\Lcal,\ H\in\Acal$ with $x(H)=-y(H)$, then there exists $z\in\Lcal$ such that $z(H)=0$ and $z(H^{\pr})=(x\circ y)(H^{\pr})$ for $H^{\pr}\in\Acal$ with $(x\circ y)(H^{\pr})=(y\circ x)(H^{\pr})$.
\end{list}

For a finite set $E$, a subset of $\{0,+,-\}^E$ satisfying (L0)-(L3) is the set of covectors of an oriented matroid.  All of the results in this paper have analogues for oriented matroids, but we stick with the language of hyperplane arrangements.

We let $\Tcal(\Acal)$ denote the set of \emph{chambers}, the maximal faces of $\Acal$.  The \emph{walls} $\Wcal(c)$ of a chamber $c$ is the set of hyperplanes in $\Acal$ incident to $c$.  Given two chambers $c,c^{\pr}\in\Tcal(\Acal)$, the \emph{separation set} $S(c,c^{\pr})$ is the set of hyperplanes in $\Acal$ separating $c$ and $c^{\pr}$; that is, $H\in S(c,c^{\pr})$ if $(c\circ c^{\pr})(H)\neq(c^{\pr}\circ c)(H)$.  Given a chamber $c_0\in\Tcal(\Acal)$, the \emph{poset of chambers} $\Pcal(\Acal,c_0)$ is an ordering on $\Tcal(\Acal)$ where $c\leq c^{\pr}$ if $S(c_0,c)\subseteq S(c_0,c^{\pr})$.  The distinguished chamber $c_0$ is called the \emph{fundamental chamber}.  If a fundamental chamber $c_0$ is given, we let $S(c)$ denote the separation set $S(c_0,c)$.  The \emph{intersection lattice} $L(\Acal)$ is the set of intersection subspaces of $\Acal$ ordered by reverse inclusion.  The set $L_k(\Acal)$ of rank $k$ elements of $L(\Acal)$ consists of all codimension $k$ intersection subspaces of $\Acal$.  The affine span $x\mapsto\bigcap\{H\in\Acal:\ x(H)=0\}$ of a face $x\in\Lcal(\Acal)$ defines an order-reversing map $\Lcal(\Acal)\ra L(\Acal)$.  For $X\in L(\Acal)$, the \emph{restriction} $\Acal^X$ is the arrangement
$$\{H\cap X:\ H\in\Acal,\ H\nsupseteq X\}$$
of hyperplanes in $X$.  If $X\in L(\Acal)$, the \emph{localization} $\Acal_X$ is the subarrangement of hyperplanes containing $X$.  For $\Acal^{\pr}\subseteq\Acal$, the restriction $x|_{\Acal^{\pr}}$ of a covector $x$ of $\Acal$ to $\Acal^{\pr}$ defines a surjective map $\Lcal(\Acal)\ra\Lcal(\Acal^{\pr})$.  Restriction also defines a surjective, order-preserving map of chamber posets $\Pcal(\Acal,c_0)\ra\Pcal(\Acal^{\pr},(c_0)|_{\Acal^{\pr}})$.  For $X\in L(\Acal)$ we let $c_X$ denote $c|_{\Acal_X}$

\begin{proposition}[see \cite{edelman:partial}]\label{prop_chamber_poset}
Let $\Acal$ be an arrangement with a fundamental chamber $c_0$.
\begin{enumerate}
\item\label{prop_chamber_poset_interval} If $X\in L(\Acal),\ x\in\Lcal(\Acal)$ with $x^{-1}(0)=\Acal_X$, then the set of chambers incident to $x$ forms an interval $[x\circ c_0,x\circ(-c_0)]$ of $\Pcal(\Acal,c_0)$ isomorphic to $\Pcal(\Acal_X,(c_0)_X)$.
\item\label{prop_chamber_poset_graded} $\Pcal(\Acal,c_0)$ is a bounded, graded poset with rank function $c\mapsto|S(c_0,c)|$.
\item\label{prop_chamber_poset_opposite} For $c,c^{\pr}\in\Tcal(\Acal)$, if $\Wcal(c)\subseteq S(c,c^{\pr})$ then $c^{\pr}=-c$.
\end{enumerate}
\end{proposition}

Combining parts \ref{prop_chamber_poset_interval} and \ref{prop_chamber_poset_opposite} of Proposition \ref{prop_chamber_poset}, we deduce the following corollary.

\begin{corollary}\label{cor_incidence_join}
Let $\Acal$ be an arrangement with fundamental chamber $c_0$ and covector $x$.  Then the join
$$\bigvee\{d\in\Pcal(\Acal,c_0)|\ d\ \mbox{covers}\ x\circ c_0\}$$
exists and is equal to $x\circ(-c_0)$.
\end{corollary}

This corollary with Lemma \ref{lem_local_lattice} implies the following.

\begin{corollary}\label{cor_local_join}
Let $\Acal$ be an arrangement with fundamental chamber $c_0$.  If $c$ is incident to $H\cap H^{\pr}$ for all chambers $c$ and hyperplanes $H,H^{\pr}\in S(c,-c_0)\cap\Wcal(c)$, then $\Pcal(\Acal,c_0)$ is a lattice.
\end{corollary}

The chambers of an arrangement completely determine its oriented matroid structure.  Mandel proved that a sign vector $x\in\{0,+,-\}^E$ is a covector of an oriented matroid with tope set $\Tcal$ if and only if $x\circ c$ is in $\Tcal$ for all $c\in\Tcal$ (see \cite[Theorem 4.2.13]{bjorner.lasVergnas.ea:oriented}).  We give a variation of this result.

\begin{theorem}\label{thm_incidence}
A chamber $c\in\Tcal(\Acal)$ is incident to $X\in L(\Acal)$ if and only if for $Y\in L(\Acal_X)$ there exists a chamber $c^{\pr}$ such that $S(c,c^{\pr})=\Acal_Y$ whenever $Y$ is incident to $c_X$.
\end{theorem}

\begin{poof}
Assume $c$ is incident to $X$, and let $x\in\Lcal(\Acal)$ such that $x^{-1}(0)=\Acal_X,\ x\leq c$.  The chamber $x\circ(-c)$ satisfies $S(c,x\circ(-c))=\Acal_X$.  Since the interval $[c,x\circ(-c)]$ is isomorphic to $\Pcal(\Acal_X,c_X)$, any wall of $c_X$ is a wall of $c$.

Now assume for $Y\in L(\Acal_X)$ there exists a chamber $c^{\pr}$ such that $S(c,c^{\pr})=\Acal_Y$ whenever $c_X$ is incident to $Y$.  We prove that $c$ is incident to $X$ by induction on the codimension of $X$.  By the inductive hypothesis, $c$ is incident to $Y$ if $Y\in L(\Acal_X)$ and $c_X$ is incident to $Y$.

Let $c^{\pr}$ be the chamber with $S(c,c^{\pr})=\Acal_X$ and let $H\in\Wcal(c^{\pr})\cap\Acal_X$.  If $H=X$, then we are done by (L3).  Thus we assume that the codimension of $X$ is at least 2.  Since $H$ is a wall of $c^{\pr}_X$, it is a wall of $c_X$.  By the inductive hypothesis, the chamber $c$ is incident to $H$.

Let $Y\in L((\Acal^H)_X),\ Y\neq H$.  If $Y\neq X$, there exists a chamber $d\in\Tcal(\Acal)$ such that $S(c,d)=\Acal_Y$ and $d$ is incident to $H$.  Hence, the chamber $d^H$ satisfies $S(c^H,d^H)=(\Acal^H)_Y$.  If $Y=X$, then $S(c^H,(c^{\pr})^H)=\Acal_X$.  By the inductive hypothesis $c^H$ is incident to $X$.  Hence, $c$ is incident to $X$.
\end{poof}

\section{Biclosed Sets}\label{sec_biclosed}

If $E$ is a finite set, a \emph{closure operator} on subsets of $E$ is a map $I\mapsto\ov{I}$ such that for $I,J\subseteq E$,
\begin{enumerate}
\item $I\subseteq\ov{I}$,
\item $\ov{\ov{I}}=\ov{I}$, and
\item $I\subseteq J$ implies $\ov{I}\subseteq\ov{J}$.
\end{enumerate}
A subset $I$ of $E$ is \emph{closed} if $I=\ov{I}$.  A closure operator is \emph{anti-exchange} if for closed subsets $I$ and $a,b\in E-I,\ a\neq b$, the element $a$ is not in $\ov{I\cup\{b\}}$ whenever $b$ is in $\ov{I\cup\{a\}}$; see \cite{edelman.jamison:convex} for many examples.  A typical example of an anti-exchange closure is the convex closure on a finite set of points in $\Rbb^n$.  We describe a disguised form of this closure operator.

\begin{definition}\label{def_convex}
Let $\Acal$ be an arrangement with fundamental chamber $c_0$, and let $I\subseteq\Acal$.
\begin{itemize}
\item $I$ is a \emph{separable set} if there exists a chamber $c\in\Tcal(\Acal)$ with $I=S(c)$.
\item $I$ is \emph{convex} if for $H\in\Acal-I$ there exists a chamber $c$ with $H\in S(c)$ and $S(c)\subseteq\Acal-I$.
\item $I$ is \emph{2-closed} if for $X\in L_2(\Acal)$ the set $I\cap\Acal_X$ is convex in $\Acal_X$ with fundamental chamber $(c_0)_X$.
\item The \emph{convex closure} (\emph{2-closure}) of $I$ is the smallest convex (2-closed) set containing $I$.
\item $I$ is \emph{biconvex} (\emph{biclosed}) if $I$ and $A\setm I$ are both convex (2-closed).
\end{itemize}
\end{definition}

Since $\Acal-S(c)=S(-c)$ holds for $c\in\Tcal(\Acal)$, separable sets are biconvex.  If $\Acal^{\pr}$ is a subarrangement of $\Acal$ and $I$ is convex in $\Acal$, then $I\cap\Acal^{\pr}$ is convex in $\Acal^{\pr}$.  We deduce the following proposition.

\begin{proposition}[see \cite{labbe:thesis}, Section 2.1]\label{prop_separable_biconvex}
If $I$ is a separable set, then $I$ is biconvex.  If $I$ is convex, then $I$ is 2-closed.  In particular, separable sets are biclosed.
\end{proposition}

If $\Acal$ is of rank 2, then separable, biconvex, and biclosed sets coincide.  Unlike the biconvex property, we prove that the biclosed property does not depend on the choice of a fundamental chamber.

\begin{lemma}\label{lem_rotation}
Let $\Acal$ be a hyperplane arrangement with chambers $c_0,c\in\Tcal(\Acal)$, and let $I\subseteq\Acal$.  The set $I$ is biclosed with respect to $c_0$ if and only if $I\vartriangle S(c_0,c)$ is biclosed with respect to $c$.
\end{lemma}

\begin{poof}
Assume $I$ is biclosed with respect to $c_0$ and let $X\in L_2(\Acal)$.  Since $I\cap\Acal_X$ is biclosed in $\Acal_X$ with respect to $(c_0)_X$, there exists a chamber $d\in\Tcal(\Acal_X)$ such that $S((c_0)_X,d)=I\cap\Acal_X$.  We have
$$S(c_X,d)=S((c_0)_X,d)\vartriangle S((c_0)_X,c_X)=(I\vartriangle S(c_0,c))\cap\Acal_X.$$
Hence, $I\vartriangle S(c_0,c)$ is biclosed with respect to $c$.
\end{poof}

\begin{lemma}\label{lem_interval_closed}
If $c,d\in\Pcal(\Acal,c_0)$ such that $c\leq d$, then $S(c,d)$ is convex.
\end{lemma}

\begin{poof}
Since $c\leq d$, the arrangement $\Acal$ has a partition into three disjoint subsets $S(c_0,c),S(c,d)$, and $S(d,-c_0)$.  Since $S(c_0,c)=S(c)$ and $S(d,-c_0)=S(-d)$, every hyperplane in $\Acal-S(c,d)$ is either in $S(c)$ or $S(d)$.
\end{poof}

Definition \ref{def_convex} has a polar dual analogue.  Given an oriented hyperplane $H$ in a real vector space $V$, let $v_H\in V^*$ be the unit vector with $v_H^{-1}(\Rbb_{>0})=H^+$.  The association $H\mapsto v_H$ defines a correspondence between (oriented) real central hyperplane arrangements and configurations of unit vectors.  Separable, convex, and 2-closed subsets of an acyclic vector configuration are usually defined as in the following proposition.

\begin{proposition}\label{prop_dictionary}
Let $I$ be a subset of hyperplanes of an arrangement $\Acal$ with fundamental chamber $c_0$.  Assume $c_0(H)=+$ for all $H\in\Acal$.
\begin{enumerate}
\item\label{prop_dictionary_separable} $I$ is a separable set if and only if there does not exist a circuit $v\in\{0,+,-\}^{\Acal}$ such that $v^{-1}(+)\subseteq I$ and $v^{-1}(-)\subseteq\Acal-I$.
\item\label{prop_dictionary_convex} $I$ is convex if and only if there does not exist a circuit $v\in\{0,+,-\}^{\Acal}$ such that $v^{-1}(+)\subseteq I,\ v^{-1}(-)\subseteq\Acal-I$ and $|v^{-1}(-)|=1$.
\item\label{prop_dictionary_2closed} $I$ is 2-closed if and only if there does not exist a circuit $v\in\{0,+,-\}^{\Acal}$ such that $v^{-1}(+)\subseteq I,\ v^{-1}(-)\subseteq\Acal-I,\ |v^{-1}(-)|=1$, and $|v^{-1}(+)|=2$.
\end{enumerate}
\end{proposition}

\begin{poof}
(\ref{prop_dictionary_separable}) For $c\in\{+,-\}^{\Acal}$, the intersection $\bigcap_{H\in\Acal}H^{c(H)}$ is nonempty if and only if $\bigcap_{H\in\Acal^{\pr}}H^{c(H)}$ is nonempty for all subarrangements $\Acal^{\pr}$ of $\Acal$.  Hence, a sign vector $c$ is a chamber if and only if there does not exist a circuit $v\in\{0,+,-\}^{\Acal}$ such that $c|_{v^{-1}(\{+,-\})}$ equals $v|_{v^{-1}(\{+,-\})}$.

(\ref{prop_dictionary_convex}) The set $I$ is convex in $\Acal$ if and only if $I$ is convex in $I\cup\{H\}$ for all $H\in\Acal-I$.  For $H\in\Acal$, $I$ is convex in $I\cup\{H\}$ if and only if there exists a chamber $c$ of $I\cup\{H\}$ such that $S(c_0|_{I\cup\{H\}},c)=\{H\}$, which holds if and only if $I$ is separable in $I\cup\{H\}$.  The statement follows from part (\ref{prop_dictionary_separable}).

(\ref{prop_dictionary_2closed}) The set $I$ is 2-closed in $\Acal$ if and only if $I\cap\Acal_X$ is convex in $\Acal_X$ for $X\in L_2(\Acal)$.  By part (\ref{prop_dictionary_convex}), this holds if and only if there does not exist a circuit $v$ of $\Acal_X$ such that $v^{-1}(+)\subseteq I,\ v^{-1}(-)\subseteq\Acal-I,\ |v^{-1}(-)=1|$ for $X\in L_2(\Acal)$.  But a circuit of $\Acal$ has three elements if and only if it is a circuit of $\Acal_X$ for some $X\in L_2(\Acal)$.
\end{poof}

\section{Galleries}\label{sec_galleries}

For chambers $c,c^{\pr}\in\Tcal(\Acal)$, a \emph{reduced gallery} from $c$ to $c^{\pr}$ is a saturated chain $c_0<\cdots<c_t$ of $\Pcal(\Acal,c)$ with $c=c_0$ and $c^{\pr}=c_t$.  A reduced gallery from $c$ to $c^{\pr}$ is \emph{incident} to a subspace $X\in L(\Acal)$ if there exists $x\in\Lcal(\Acal),\ x^{-1}(0)=\Acal_X$ such that the chambers $x\circ c$ and $x\circ(-c)$ are both in the gallery.  If $r=c_0,\ldots,c_t$ is a reduced gallery, we let $r_X$ denote the reduced gallery $(c_0)_X,\ldots,(c_t)_X$ in $\Acal_X$.  For reduced galleries $r,r^{\pr}$ from $c$ to $c^{\pr}$, the $L_2$-separation set $L_2(r,r^{\pr})$ is the set of codimension 2 subspaces for which $r_X\neq r^{\pr}_X$.

\begin{proposition}\label{prop_gallery_incidence}
If $r$ is a reduced gallery from $c$ to $c^{\pr}$ and $X$ a codimension 2 intersection subspace, then there exists a reduced gallery $r^{\pr}$ such that $L_2(r,r^{\pr})=\{X\}$ if and only if $r$ is incident to $X$.
\end{proposition}

\begin{poof}
Assume $r$ is incident to $X$.  Let $x\in\Lcal(\Acal)$ with $x^{-1}(0)=\Acal_X$ such that $x\circ c$ and $x\circ(-c)$ are chambers in $r$.  By Proposition \ref{prop_chamber_poset}(\ref{prop_chamber_poset_interval}), the interval $[x\circ c,x\circ(-c)]$ of $\Pcal(\Acal,c)$ is isomorphic to $\Pcal(\Acal_X,c_X)$.  Since $\Acal_X$ is of rank 2, the interval $[x\circ c,x\circ(-c)]$ has two maximal chains.

Let $r^{\pr}$ be the symmetric difference of $r$ with the open interval $(x\circ c,x\circ(-c))$.  Then $r^{\pr}$ is a reduced gallery from $c$ to $c^{\pr}$ such that $X\in L_2(r,r^{\pr})$.  Let $Y\in L_2(\Acal),\ Y\neq X$.  Since every hyperplane in $S(x\circ c,x\circ(-c))$ contains $X$, the localized arrangement $\Acal_Y$ contains at most one hyperplane of $S(x\circ c,x\circ(-c))$.  Hence, $d_Y=(x\circ c)_Y$ or $d_Y=x\circ(-c)_Y$ for $d\in[x\circ c,x\circ(-c)]$, so $L_2(r,r^{\pr})=\{X\}$.

Now assume $r^{\pr}$ is a reduced gallery from $c$ to $c^{\pr}$ such that $L_2(r,r^{\pr})=\{X\}$.  Let $d$ be the largest chamber common to $r$ and $r^{\pr}$ for which $r_{\leq d}=r^{\pr}_{\leq d}$.  Let $H$ and $H^{\pr}$ be the upper walls of $d$ crossed by $r$ and $r^{\pr}$.  Since $r$ and $r^{\pr}$ are separated by $H\cap H^{\pr}$, both $H$ and $H^{\pr}$ contain $X$.  Let $x\in\{0,+,-\}^{\Acal}$ such that $x(H^{\prpr})=0$ if $H^{\prpr}\in\Acal_X$ and $x(H^{\prpr})=d(H^{\prpr})$ otherwise.

If $x\circ(-c)$ is not a chamber in $r$, then there exists a hyperplane $H^{\prpr}\in\Acal$ not containing $X$ such that $r$ crosses $H^{\prpr}$ before $H^{\pr}$ but after $H$.  Then $r$ and $r^{\pr}$ are separated by $H^{\pr}\cap H^{\prpr}$, an impossibility.  Hence, $x\circ(-c)$ is a chamber.  By Theorem \ref{thm_incidence}, we conclude that $x$ is a covector of $\Acal$ and $r$ is incident to $x$.
\end{poof}

\begin{figure}
\begin{centering}
\begin{tikzpicture}[scale=1.2]
\begin{scope}[yshift=2.5cm]
  \draw (0,0) circle (1.8cm);
  \draw[rotate=60] (1.8cm,0) arc (0:180:1.8cm and .8cm);
  \draw[rotate=180] (1.8cm,0) arc (0:180:1.8cm and .8cm);
  \draw[rotate=300] (1.8cm,0) arc (0:180:1.8cm and .8cm);
  \draw (270:1.4cm) node{$c_0$};
  \draw (270:1.8cm) node{$A$}
        (270:.8cm) node{$D$}
        (150:.8cm) node{$B$}
        (30:.8cm) node{$C$};
\end{scope}

\begin{scope}[xshift=5cm,yshift=.5cm]
  \coordinate (O) at (0,0);
  \coordinate (A) at (-1.5,0.8);
  \coordinate (B) at (-0.5,1.2);
  \coordinate (C) at (0.5,1.2);
  \coordinate (D) at (1.5,0.8);
  \draw (O) -- node{$A$} (A);
  \draw (O) -- node{$B$} (B);
  \draw (O) -- node{$D$} (C);
  \draw (O) -- node{$C$} (D);
  \draw (A) -- ($(A)+(B)$);
  \draw (B) -- ($(A)+(B)$);
  \draw (B) -- ($(B)+(C)$);
  \draw (C) -- ($(B)+(C)$);
  \draw (C) -- ($(C)+(D)$);
  \draw (D) -- ($(C)+(D)$);
  \draw ($(A)+(B)$) -- ($(A)+(B)+(C)$);
  \draw ($(B)+(C)$) -- ($(A)+(B)+(C)$);
  \draw ($(B)+(C)$) -- ($(B)+(C)+(D)$);
  \draw ($(C)+(D)$) -- ($(B)+(C)+(D)$);
  \draw ($(A)+(B)+(C)$) -- ($(A)+(B)+(C)+(D)$);
  \draw ($(B)+(C)+(D)$) -- ($(A)+(B)+(C)+(D)$);
  \begin{scope}[black!70]
    \draw (A) -- ($(A)+(D)$);
    \draw (D) -- ($(A)+(D)$);
    \draw ($(A)+(B)$) -- ($(A)+(B)+(D)$);
    \draw ($(C)+(D)$) -- ($(A)+(C)+(D)$);
    \draw ($(A)+(D)$) -- ($(A)+(B)+(D)$);
    \draw ($(A)+(D)$) -- ($(A)+(C)+(D)$);
    \draw ($(A)+(B)+(D)$) -- ($(A)+(B)+(C)+(D)$);
    \draw ($(A)+(C)+(D)$) -- ($(A)+(B)+(C)+(D)$);
  \end{scope}
\end{scope}

\begin{scope}[xshift=10cm, yscale=.8]
  \draw[shorten >=5pt, shorten <=5pt] (0,0) -- (-1,1);
  \draw[shorten >=2pt, shorten <=2pt] (0,0) -- (0,1);
  \draw[shorten >=5pt, shorten <=5pt] (0,0) -- (1,1);
  \draw[shorten >=2pt, shorten <=2pt] (-1,1) -- (-1,2);
  \draw[shorten >=5pt, shorten <=5pt] (-1,1) -- (0,2);
  \draw[shorten >=5pt, shorten <=5pt] (0,1) -- (1,2);
  \draw[shorten >=2pt, shorten <=2pt] (1,1) -- (1,2);
  \draw[shorten >=3pt, shorten <=3pt] (-1,2) -- (-.5,3);
  \draw[shorten >=3pt, shorten <=3pt] (0,2) -- (-.5,3);
  \draw[shorten >=3pt, shorten <=3pt] (1,2) -- (.5,3);
  \draw[shorten >=3pt, shorten <=3pt] (-.5,3) -- (-1,4);
  \draw[shorten >=3pt, shorten <=3pt] (.5,3) -- (0,4);
  \draw[shorten >=3pt, shorten <=3pt] (.5,3) -- (1,4);
  \draw[shorten >=2pt, shorten <=2pt] (-1,4) -- (-1,5);
  \draw[shorten >=5pt, shorten <=5pt] (-1,4) -- (0,5);
  \draw[shorten >=5pt, shorten <=5pt] (0,4) -- (1,5);
  \draw[shorten >=2pt, shorten <=2pt] (1,4) -- (1,5);
  \draw[shorten >=5pt, shorten <=5pt] (-1,5) -- (0,6);
  \draw[shorten >=2pt, shorten <=2pt] (0,5) -- (0,6);
  \draw[shorten >=5pt, shorten <=5pt] (1,5) -- (0,6);
  \draw (0,0) node{$\ss ABCD$}
        (-1,1) node{$\ss ACBD$}
        (0,1) node{$\ss BACD$}
        (1,1) node{$\ss ABDC$}
        (-1,2) node{$\ss CABD$}
        (0,2) node{$\ss ACDB$}
        (1,2) node{$\ss BADC$}
        (-.5,3) node{$\ss CADB$}
        (.5,3) node{$\ss BDAC$}
        (-1,4) node{$\ss CDAB$}
        (0,4) node{$\ss BDCA$}
        (1,4) node{$\ss DBAC$}
        (-1,5) node{$\ss CDBA$}
        (0,5) node{$\ss DCAB$}
        (1,5) node{$\ss DBCA$}
        (0,6) node{$\ss DCBA$};
\end{scope}

\end{tikzpicture}
\caption{\label{fig_gallery}\scriptsize(left) An arrangement of four planes in $\Rbb^3$. (center) The poset of chambers. (right) The graph of reduced galleries from $c_0$ to $-c_0$.  The galleries $BACD$ and $CDBA$ are separated by four codimension 2 subspaces, but the shortest path between them in the gallery graph has length six.}
\end{centering}
\end{figure}

For $c,c^{\pr}\in\Tcal(\Acal)$, the set of reduced galleries from $c$ to $c^{\pr}$ forms a graph $G_2(c,c^{\pr})$ where galleries $r$ and $r^{\pr}$ are adjacent if $|L_2(r,r^{\pr})|=1$.  The function $L_2(\cdot,\cdot):G_2\times G_2\ra\Zbb$ satisfies
\begin{list}{}{}
\item[(M1)] $L_2(r,r^{\pr})=L_2(r^{\pr},r)$,
\item[(M2)] $r$ is adjacent to $r^{\pr}$ implies $|L_2(r,r^{\pr})|=1$, and
\item[(M3)] $L_2(r,r^{\pr})=L_2(r,r^{\prpr})\vartriangle L_2(r^{\prpr},r^{\pr})$,
\end{list}
for reduced galleries $r,r^{\pr},r^{\prpr}$.  Here, $X\vartriangle Y$ is the symmetric difference $(X\setminus Y)\cup(Y\setminus X)$ for sets $X,Y$.  A function satisfying (M1),(M2),(M3) on some graph $G$ is called a \emph{set-valued metric} on $G$.  The function $S:\Tcal(\Acal)\times\Tcal(\Acal)\ra\Zbb$ is also a set-valued metric.

The separation function $S(\cdot,\cdot)$ has the additional property that for $c,c^{\pr}\in\Tcal(\Acal)$ there exists a geodesic from $c$ to $c^{\pr}$ in the chamber graph of length $|S(c,c^{\pr})|$.  However, this property does not hold for $L_2$; see Figure \ref{fig_gallery}.  In general, the length of the shortest geodesic between two galleries is bounded below by the size of their $L_2$-separation set.  Still, an analogue of Proposition \ref{prop_chamber_poset}(\ref{prop_chamber_poset_graded}) holds for some gallery graphs, as described in Proposition \ref{prop_accessible}.  A reduced gallery $r_0$ from $c$ to $c^{\pr}$ is \emph{$L_2$-accessible} if there exists a geodesic between $r_0$ and $r$ of length $|L_2(r_0,r)|$ for all $r\in G_2(c,c^{\pr})$.

\begin{proposition}[\cite{reiner.roichman:diameter}, Proposition 3.12]\label{prop_accessible}
If $G_2(c,-c)$ contains an accessible gallery, then the diameter of $G_2(c,-c)$ is equal to $|L_2(\Acal)|$.
\end{proposition}

A permutation $\pi:H_1,H_2,\ldots$ of $\Acal$ is \emph{admissible} if for each codimension 2 subspace $X$, there exists a gallery of $\Acal_X$ from $(c_0)_X$ to $-(c_0)_X$ crossing the hyperplanes in the order defined by $\pi$.  If $c_0,c_1,\ldots$ is a gallery of $\Acal$ then $H_1,H_2,\ldots$ is an admissible permutation of $\Acal$ where $S(c_{i-1},c_i)=\{H_i\}$.  The following proposition gives a partial converse.

\begin{proposition}\label{prop_admissible}
Let $\Acal$ be an arrangement with fundamental chamber $c_0$ such that every biclosed set is the separation set of some chamber.  If $\pi:H_1,H_2,\ldots,H_N$ is an admissible permutation of $\Acal$ then there exists a gallery $c_0,c_1,\ldots,c_N=-c_0$ such that $S(c_{i-1},c_i)=\{H_i\}$ for all $i$.
\end{proposition}

\begin{poof}
For $0\leq j\leq N$, let $I_j=\{H_1,\ldots,H_j\}$.  Let $X\in L_2(\Acal)$, and suppose
$$\Acal_X=\{H_{p_1},\ldots,H_{p_t}\},\hspace{5mm} p_1<\cdots<p_t.$$
Since $\pi$ is admissible, there exists a chamber $c$ of $\Acal_X$ such that $S((c_0)_X,c)=\{H_{p_k}\ |\ p_k\leq i\}=I_j\cap\Acal_X$.  Hence, $I_j$ is biclosed, and there exists a chamber $c_j$ such that $S(c_j)=I_j$.  The chain $c_0<c_1<\cdots<c_N$ is a reduced gallery of $\Acal$ inducing $\pi$.
\end{poof}

For the generic arrangement of four planes in $\Rbb^3$ shown in Figure \ref{fig_gallery}, all 16 subsets of hyperplanes are biclosed, but there are only 14 chambers.  This example shows the necessity of the conditions in Proposition \ref{prop_admissible}, as all 24 permutations of the hyperplanes are admissible, but only 16 come from reduced galleries.

Admissible permutations are easier to flip than reduced galleries, as described in the following lemma.

\begin{lemma}\label{lem_admissible_flip}
Let $H_1,H_2,\ldots,H_N$ be an admissible permutation of $\Acal$.  Suppose the set of hyperplanes containing some codimension 2 subspace $X$ is a contiguous subsequence $H_i,H_{i+1},\ldots,H_j$.  Then the permutation
$$H_1,\ldots,H_{i-1},H_j,\ldots,H_i,H_{j+1},\ldots,H_N$$
obtained by flipping the subsequence $H_i,\ldots,H_j$ is an admissible permutation.
\end{lemma}

\begin{poof}
Let $\pi:H_1,\ldots,H_N$ be the original permutation and let $\pi^{\pr}$ be the flip.  If $Y$ is some codimension 2 subspace not contained in at least two hyperplanes in $H_i,\ldots,H_j$, then the restrictions of $\pi$ and $\pi^{\pr}$ to $\Acal_X$ are the same.  If $Y$ does contain at least two hyperplanes in the sequence, then it equals $X$.  Since $H_i,\ldots,H_j$ is an admissible permutation of $\Acal_X$, so is the reverse $H_j,\ldots,H_i$.
\end{poof}

Once again, the arrangement of Figure \ref{fig_gallery} gives an example of a gallery that cannot be flipped as in Lemma \ref{lem_admissible_flip}.  For instance, the gallery from $c_0$ to $-c_0$ crossing hyperplanes $B,D,C,A$ in that order is not incident to $C\cap D$, so it does not admit a flip to $B,C,D,A$.  However, the permutation $B,C,D,A$ is an admissible permutation of $\Acal$.

Theorem \ref{thm_main_gallery} now follows immediately from Proposition \ref{prop_admissible} and Lemma \ref{lem_admissible_flip}.

\section{Simplicial Arrangements}\label{sec_simplicial}

An arrangement is \emph{simplicial} if every chamber is a simplicial cone.  An arrangement $\Acal$ with fundamental chamber $c_0$ is \emph{bineighborly} if for $c\in\Pcal(\Acal,c_0),\ H,H^{\pr}\in\Wcal(c)\cap S(c)$, the chamber $c$ is incident to $H\cap H^{\pr}$.  In separate work, we proved that $(\Acal,c_0)$ is bineighborly if and only if $\Pcal(\Acal,c_0)$ is a semidistributive lattice \cite{mcconville:crosscut}.

Our proof of Theorem \ref{thm_bineighborly} is similar to Dyer's proof in the root system case (see \cite[Theorem 1.5]{dyer:weak}).

\begin{theorem}\label{thm_bineighborly}
If $\Acal$ is a bineighborly arrangement with respect to a fundamental chamber $c_0$, then $S(x\vee y)$ equals the 2-closure of $S(x)\cup S(y)$ for $x,y\in\Pcal(\Acal,c_0)$.
\end{theorem}

\begin{poof}
By Lemma \ref{lem_interval_closed}, it suffices to prove the following claim.

Claim: If $a\in\Pcal(\Acal,c_0)$ and $x,y\in[a,-c_0]$, then $S(a,x\vee y)\subseteq\ov{S(a,x)\cup S(a,y)}$ holds.

If $a=-c_0$, then the claim is immediate.

Let $a\in\Pcal(\Acal,c_0),\ a<-c_0$ and assume the claim holds if $a$ is replaced by $a^{\pr}\in\Pcal(\Acal,c_0)$ with $a<a^{\pr}$.  Let $x,y\in[a,-c_0]$.  We may assume $x\wedge y=a$ as otherwise we have
$$S(a,x\vee y)=S(a,x\wedge y)\cup S(x\wedge y)=S(a,x\wedge y)\cup\ov{S(x\wedge y,x)\cup S(x\wedge y,y)}\subseteq\ov{S(a,x)\cup S(a,y)}.$$

If $x=a$ then the identity
$$S(a,x\vee y)=S(a,y)=\ov{S(a,y)}=\ov{S(a,x)\cup S(a,y)}$$
is clear.

Assume $a<x$ and $a<y$ hold.  Let $H\in U(a)\cap S(a,x)$ and $H^{\pr}\in U(a)\cap S(a,y)$.  Let $c,d$ denote the chambers covering $a$ with $S(a,c)=\{H\}$ and $S(a,d)=\{H^{\pr}\}$.  Since $(\Acal,c_0)$ is bineighborly, $a$ is incident to $H\cap H^{\pr}$, so the join $c\vee d$ satisfies
$$S(a,c\vee d)=\Acal_{H\cap H^{\pr}}=\ov{S(a,c)\cup S(a,d)}.$$

The equality 
$$S(c\vee d,x\vee y)=\ov{S(c\vee d,x\vee d)\cup S(c\vee d,c\vee y)}$$
holds by the induction hypothesis.  The rest of this string of equalities and inequalities follows from properties of closure operators.

\begin{align*}
\ov{S(a,x)\cup S(a,y)}&=\ov{S(a,c)\cup S(c,x)\cup S(a,d)\cup S(d,y)}\\
&=\ov{S(a,c\vee d)\cup S(c,x)\cup S(d,y)}\\
&=\ov{S(c,x)\cup S(c,c\vee d)\cup S(d,y)\cup S(d,c\vee d)}\\
&=\ov{S(c,x\vee d)\cup S(d,c\vee y)}\\
&\supseteq S(a,c\vee d)\cup\ov{S(c\vee d,x\vee d)\cup S(c\vee d,c\vee y)}\\
&=S(a,c\vee d)\cup S(c\vee d,x\vee y)\\
&=S(a,x\vee y)
\end{align*}

\end{poof}

Reading defined an arrangement to be \emph{bisimplicial} if $c|_{\Wcal(c)\cap S(c)}$ is a simplicial cone for all $c\in\Tcal(\Acal)$ \cite{reading:lattice_congruence}.  We prove in Proposition \ref{prop_bisimplicial} bineighborly arrangements are bisimplicial.  This result is somewhat surprising since there exist non-simplicial neighborly polytopes.  Not all bisimplicial arrangements define semidistributive lattices nor have joins computed by a 2-closure; see Figure \ref{fig_gallery}.

\begin{proposition}\label{prop_bisimplicial}
Let $\Acal$ be an arrangement with fundamental chamber $c_0$.  If $\Acal$ is bineighborly, then $\Acal$ is bisimplicial.
\end{proposition}

\begin{poof}
Let $c$ be a chamber of $\Acal$, and let $\Acal^{\pr}=\Wcal(c)\cap S(c,-c_0)$.  By the proof of Theorem \ref{thm_bineighborly}, for $I\subseteq\Acal^{\pr}$, the set of hyperplanes separating $c$ and
$$\bigvee\{c^{\pr}:\ (\exists H\in I)\ S(c,c^{\pr})=\{H\}\}$$
is equal to the 2-closure of $I$ in $\Acal$.  Since every subset of $\Acal^{\pr}$ is 2-closed in $\Acal^{\pr}$, this implies $\Acal^{\pr}$ has $2^{|\Acal^{\pr}|}$ distinct chambers, so $c|_{\Acal^{\pr}}$ is simplicial.
\end{poof}

For $H\in\Acal$, let $\depth(H)$ be the minimum size of $S(c)$ where $H\in S(c),\ c\in\Pcal(\Acal,c_0)$.

\begin{proposition}\label{prop_depth_simplicial}
Let $\Acal$ be a simplicial arrangement with fundamental chamber $c_0$, and let $I\subseteq\Acal$.  If $I$ is 2-closed and $I\supseteq\Wcal(c_0)$, then $I=\Acal$.
\end{proposition}

\begin{poof}
If $\depth(H)=1$, then $H\in I$ by assumption.

Let $k>1$ and suppose $H\in I$ if $\depth(H)<k$.  Let $H\in\Acal$ such that $\depth(H)=k$.  Let $c$ be a chamber such that $|S(c)|=k$ and $H\in S(c)$.  By minimality of $c$, $H$ is the unique wall of $c$ separating it from $c_0$.  Let $d$ be the chamber with $S(d,c)=\{H\}$, and let $H^{\pr}$ be any hyperplane in $\Wcal(d)\cap S(d)$.  Since $d$ is simplicial, it is incident to $H\cap H^{\pr}$.

Let $x$ be the codimension 2 face supported by $H\cap H^{\pr}$ incident to $d$.  Let $H_1,H_2$ be the two walls of $x\circ c_0$ containing $x$.  Let $c_1,c_2$ be the chambers with $S(x\circ c_0,c_i)=\{H_i\}$ for $i=1,2$.  Since $|S(c_i)|<k$ for $i=1,2$, the depths of $H_1$ and $H_2$ are both less than $k$.  Hence, $H_1,H_2\in I$.  Since $I$ is 2-closed, this implies $H\in I$.
\end{poof}

\begin{theorem}
Let $\Acal$ be a simplicial arrangement with fundamental chamber $c_0$.  A subset $I$ of $\Acal$ is biclosed if and only if $I=S(c)$ for some chamber $c$.
\end{theorem}

\begin{poof}
Let $I\subseteq\Acal$.  If $I=S(c)$ for some chamber $c$, then $I$ is biclosed by Proposition \ref{prop_separable_biconvex}.

Assume $I$ is biclosed.  Choose a chamber $c$ minimizing $|I\vartriangle S(c)|$.  Let $I^{\pr}=I\vartriangle S(c)$.  By Lemma \ref{lem_rotation} $I^{\pr}$ is biclosed with respect to $c$.  The minimality of $c$ implies $\Wcal(c)\subseteq\Acal-I^{\pr}$.  Since $\Acal-I^{\pr}$ is 2-closed with respect to $c$ and it contains $\Wcal(c)$, Proposition \ref{prop_depth_simplicial} implies $I^{\pr}$ is empty.
\end{poof}

\section{Supersolvable Arrangements}\label{sec_supersolvable}

An intersection subspace $X\in L(\Acal)$ of an arrangement $\Acal$ is \emph{modular} if $X+Y$ is in $L(\Acal)$ for all $Y\in L(\Acal)$.  A rank $r$ arrangement $\Acal$ is \emph{supersolvable} if its intersection lattice contains a maximal chain $X_0<X_1<\cdots<X_r$ where $X_i\in L_i(\Acal)$ is modular.  If $X\in L_{r-1}(\Acal)$ for a rank $r$ arrangement $\Acal$, then $X$ is called a \emph{modular line}.

Most results about supersolvable arrangements are proved inductively by localization at a modular line.  This approach to supersolvable arrangements is suggested by the following recursive characterization obtained by Bj{\"o}rner, Edelman, and Ziegler.

\begin{theorem}\label{thm_super_recur}(BEZ \cite{bjorner.edelman.ziegler:lattice})
Every arrangement of rank at most 2 is supersolvable.  An arrangement $\Acal$ of rank $r\geq 3$ is supersolvable if and only if it contains a modular line $l$ such that the localization $\Acal_l$ is supersolvable.
\end{theorem}

An extension of an arrangement $\Acal$ is a hyperplane arrangement containing $\Acal$.  Although supersolvable arrangements are far from generic, any arrangement can be made supersolvable by adding enough hyperplanes.

\begin{corollary}\label{cor_super_recur}
Any arrangement admits a supersolvable extension of the same rank.
\end{corollary}

\begin{poof}
Let $\Acal$ be an arbitrary arrangement of rank $r\geq 3$, and let $l\in L_{r-1}(\Acal)$ be a line.  Decompose $\Acal$ into a disjoint union $\Acal=\Acal_l\sqcup(\Acal\setm\Acal_l)$.  Let $\Acal_0$ be the union of $\Acal_l$ with the set of hyperplanes of the form $l+(H\cap H^{\pr})$ for some pair of hyperplanes $H,H^{\pr}$ in $\Acal\setm\Acal_l$.  By the inductive hypothesis, $\Acal_0$ has a supersolvable extension $\wtil{\Acal_0}$ of rank $r-1$.  Then the disjoint union $\wtil{\Acal_0}\sqcup(\Acal\setm\Acal_l)$ is a rank $r$ supersolvable arrangement by Theorem \ref{thm_super_recur}.
\end{poof}

The following proposition due to Reiner and Roichman was essential to computing the diameter of the graph of reduced galleries in a supersolvable arrangement.

\begin{proposition}[\cite{reiner.roichman:diameter}, Proposition 4.6]\label{prop_fiber_cham}
Assume that $l$ is a modular line of an arrangement $\Acal$ with chamber $c_0$ incident to $l$.  Let $\pi:\Pcal(\Acal,c_0)\ra\Pcal(\Acal_l,(c_0)_l)$ be the localization map, and let $U$ denote the fiber $\pi^{-1}(\pi(c_0))$.
\begin{enumerate}
\item\label{prop_fiber_cham_p2} The fiber $U=\{c_0,c_1,\ldots,c_t\}$ is linearly ordered $c_0<c_1<\cdots<c_t$.  This induces a linear order $H_1,H_2,\ldots,H_t$ on $\Acal\setm\Acal_l$ such that $H_i$ is the unique hyperplane in $S(c_{i-1},c_i)$.
\item\label{prop_fiber_cham_p3} Using the linear order on $\Acal\setm\Acal_l$ from part (\ref{prop_fiber_cham_p2}), if $i<j<k$ and if the chamber $c_0$ incident to $l$ is also incident to $l+H_i\cap H_k$ then $H_j\supseteq H_i\cap H_k$.
\end{enumerate}
\end{proposition}

\begin{theorem}\label{thm_main_super}
Let $\Acal$ be a supersolvable arrangement with fundamental chamber $c_0$ incident to a modular flag.
\begin{enumerate}
\item For $I\subseteq\Acal$, there exists a chamber $c\in\Tcal(\Acal)$ with $I=S(c_0,c)$ if and only if $I$ is biclosed.
\item For $c,d\in\Tcal(\Acal)$, the separation set $S(c\vee d)$ is the 2-closure of $S(c)\cup S(d)$.
\end{enumerate}
\end{theorem}

\begin{poof}
Let $l\in L(\Acal)$ be a modular line incident to $c_0$ such that $(c_0)_l$ is incident to a modular flag of $\Acal_l$.  Assume both parts of the theorem hold for the pair $(\Acal_l,(c_0)_l)$.

(1) Let $I\subseteq\Acal$.  If $I=S(c_0,c)$ for some chamber $c$, then $I$ is biclosed by Proposition \ref{prop_separable_biconvex}.

Assume $I$ is biclosed.  The restriction $I\cap\Acal_l$ is $(c_0)_l$-biclosed, so there exists a chamber $\ov{c}\in\Tcal(\Acal_l)$ such that $S((c_0)_l,\ov{c})=I\cap\Acal_l$.  Since $c_0$ is incident to $l$, there exists a chamber $c\in\Tcal(\Acal)$ such that $S(c_0,c)=I\cap\Acal_l$ by Proposition \ref{prop_chamber_poset}(\ref{prop_chamber_poset_interval}).  Let $c_1,\ldots,c_{t+1}\in\Tcal(\Acal)$ such that $c=c_1,\ S(c_1,c_{t+1})=\Acal-\Acal_l,\ (c_i)_l=\ov{c}$ and $|S(c_i,c_{i+1})|=1$ for all $i$.  Let $H_1,\ldots,H_t$ be the hyperplanes of $\Acal-\Acal_l$ where $S(c_i,c_{i+1})=\{H_i\}$ for all $i$.

Assume $H_i\in I$ for some $i>1$, and let $X=H_{i-1}\cap H_i$.  Since $l$ is modular, $X+l$ is a hyperplane of $\Acal$ containing $l$.  If $X+l\in S(c_0,c)$ then $H_{i-1}$ is in the 2-closure of $\{X+l,H_i\}\subseteq I$.  If $X+l\notin S(c_0,c)$ then $H_i$ is in the 2-closure of $\{X+l,H_{i-1}\}$.  Since $I$ is biclosed, both cases imply $H_{i-1}\in I$.  Hence, $I\cap(\Acal-\Acal_l)$ is an initial segment of hyperplanes $H_1,\ldots,H_k$ for some $k$, so $I\cap(\Acal-\Acal_l)=S(c,c_{k+1})$ holds.  Therefore, we obtain
$$I=S(c_0,c)\cup S(c,c_{k+1})=S(c_0,c_{k+1}).$$

(2) Let $c,d\in\Tcal(\Acal)$.  The equality
$$S((c\vee d)_l)=\ov{S(c_l)\cup S(d_l)}$$
holds by the assumption on $(\Acal_l,(c_0)_l)$.  Since $c_0$ is incident to $l$, there exist chambers $b,c^{\pr},d^{\pr}\in\Tcal(\Acal)$ incident to $l$ such that $b_l=(c\vee d)_l,\ c^{\pr}_l=c_l,\ d^{\pr}_l=d_l$ by Proposition \ref{prop_chamber_poset}(\ref{prop_chamber_poset_interval}).  The above equality then lifts to
$$S(b)=\ov{S(c^{\pr})\cup S(d^{\pr})}.$$

Let $c_1,\ldots,c_{t+1}\in\Tcal(\Acal)$ such that $c=c_1,\ S(c_1,c_{t+1})=\Acal-\Acal_l,\ (c_i)_l=\ov{b}$ and $|S(c_i,c_{i+1})|=1$ for all $i$.  Let $H_1,\ldots,H_t$ be the hyperplanes of $\Acal-\Acal_l$ where $S(c_i,c_{i+1})=\{H_i\}$ for all $i$.

The join $c\vee d$ is equal to $c_{k+1}$ for some $k$.  Since $H_k$ is a wall of $c\vee d$, either $H_k\in S(c)$ or $H_k\in S(d)$.  By symmetry, we may assume $H_k\in S(c)$.  Suppose there exists $H_i$ with $i<k$, and let $X=H_i\cap H_k,\ H=X+l$.  If $H\in S(b)$ then $H_i$ is in the 2-closure of $\{H,H_k\}$, so it lies in $\ov{S(c)\cup S(d)}$.  If $H\notin S(b)$, then $H\notin S(c)$ and $H_k$ is in the 2-closure of $\{H_i,H\}$.  Since $S(c)$ is biclosed, this forces $H_i\in S(c)$.
\end{poof}

The following corollary follows immediately from Theorem \ref{thm_main_super} and Proposition \ref{prop_admissible}.

\begin{corollary}\label{cor_admiss_gal}
Let $\Acal$ be a supersolvable arrangement with fundamental chamber $c_0$.  For any admissible permutation $H_1,\ldots,H_N$ of $\Acal$ there exists a reduced gallery $c_0,\ldots,c_N$ such that $S(c_{i-1},c_i)=\{H_i\}$.
\end{corollary}

We now fill in the gap in the proof of \cite[Theorem 1.1]{reiner.roichman:diameter} by adapting that proof in our language, and applying Corollary \ref{cor_admiss_gal} with Lemma \ref{lem_admissible_flip}.

\begin{theorem}[\cite{reiner.roichman:diameter}, Theorem 1.1]\label{thm_diameter}
Let $\Acal$ be a rank $n$ supersolvable arrangement with fundamental chamber $c_0$ incident to a modular flag.  The graph of reduced galleries from $c_0$ to $-c_0$ has diameter $|L_2(\Acal)|$.
\end{theorem}

\begin{poof}
Let $X_1\supsetneq\cdots\supsetneq X_{n-1}$ be a modular flag such that $c_0$ is incident to $X_i\in L_i(\Acal)$ for all $i$.  Let $x_1,\ldots,x_{n-1}$ be the unique covectors of $\Acal$ with $c_0\geq x_1\geq\cdots\geq x_{n-1}$ and $x_i^{-1}(0)=\Acal_{X_i}$ for all $i$.  We set $l=X_{n-1}$.  Let $r_0$ be a maximal chain in $\Pcal(\Acal,c_0)$ extending the chain $x_1\circ(-c_0)<\cdots<x_{n-1}\circ(-c_0)$.  By Proposition \ref{prop_accessible}, it suffices to show that $r_0$ is $L_2$-accessible.  We proceed by induction on $n$.

Let $r$ be a reduced gallery from $c_0$ to $-c_0$ and assume $L_2(r_0,r)\subseteq L_2(\Acal_l)$.  If $H\in\Acal_l,\ H^{\pr}\in\Acal-\Acal_l$ then $(r_0)_{H\cap H^{\pr}}=r_{H\cap H^{\pr}}$ and $(r_0)_{H\cap H^{\pr}}$ crosses $H$ before $H^{\pr}$.  Hence, $r$ must contain the chamber $x_{n-1}\circ(-c_0)$.  By Proposition \ref{prop_fiber_cham}(\ref{prop_fiber_cham_p2}), the interval $[x_{n-1}\circ(-c_0),-c_0]$ of $\Pcal(\Acal,c_0)$ is a chain, so the galleries $r$ and $r_0$ agree above $x_{n-1}$.  Since $[c_0,x_{n-1}\circ(-c_0)]$ is isomorphic to $\Pcal(\Acal_l,(c_0)_l)$, the distance between $r_0$ and $r$ in the reduced gallery graph is equal to $|L_2(r_0,r)|$ by the induction hypothesis.

Now assume $L_2(r_0,r)\nsubseteq L_2(\Acal_l)$.  Let $H_1,\ldots,H_N$ be the total order on $\Acal$ induced by $r$, and let $k$ be the smallest index for which $H_k\supseteq l$ but $H_{k-1}\nsupseteq l$.  Let $X=H_{k-1}\cap H_k$.  Let $c_1$ be the largest chamber in $r$ incident to $x_{n-1}$.  By the assumption on $k$, $(c_1)_l$ is incident to $H_k$.  Since $c_1$ is incident to $l$, this implies $c_1$ is incident to $H_k$.  If $i,j$ are indices with $i<j<k$ such that $H_i\cap H_k=X$ holds, then $H_j\supseteq X$ by Proposition \ref{prop_fiber_cham}(\ref{prop_fiber_cham_p3}).  By Corollary \ref{cor_admiss_gal}, the gallery $r$ is incident to $X$, so there exists a gallery $r^{\pr}$ adjacent to $r$ such that $|L_2(r_0,r^{\pr})|=|L_2(r_0,r)|-1$.  By induction, there exists a gallery $r^{\prpr}$ of distance $|L_2(r_0,r)|-|L_2((r_0)_l,r_l)|$ from $r$ such that $L_2(r_0,r^{\prpr})\subseteq L_2(\Acal_l)$.  By the previous case, this implies the distance between $r_0$ and $r$ is equal to $|L_2(r_0,r)|$.
\end{poof}

\section{Acknowledgements}

The author thanks Josh Hallam for helpful discussions on supersolvable lattices.  He thanks Jean-Philippe Labb{\'e} for suggesting the reference \cite{pilkington:convex}, and he thanks Christophe Hohlweg and Jean-Philippe Labb{\'e} for describing their perspectives on biclosed sets in infinite root systems.  He also thanks Vic Reiner and Pasha Pylyavskyy for their help and guidance.

\scriptsize
\bibliography{bib_biclosed}{}
\bibliographystyle{plain}

\end{document}